\documentclass[12pt]{article}
\usepackage{amsfonts}
\usepackage{amsmath, amssymb, amsthm}
\usepackage{enumerate}
\usepackage{epsfig}
\usepackage{verbatim}

\input{tcilatex}
\begin{document}

\title{Isogeometric Analysis for singularly perturbed problems
in 1-D: a numerical study\\
}
\author{K. Liotati and C. Xenophontos\thanks{
Corresponding author. Email: xenophontos@ucy.ac.cy} \\
Department of Mathematics and Statistics \\
University of Cyprus \\
P.O. Box 20537 \\
1678 Nicosia \\
Cyprus}
\maketitle

\begin{abstract}
We perform numerical experiments on one-dimensional singularly perturbed
problems of reaction-convection-diffusion type, using isogeometric analysis.
In particular, we use a Galerkin formulation with B-splines as basis
functions. The question we address is: how should the knots be chosen in
order to get uniform, exponential convergence in the maximum norm? We
provide specific guidelines on how to achieve precisely this, for three
different singularly perturbed problems.
\end{abstract}

\section{Introduction}
\label{intro}

We consider second order singularly perturbed problems (SPPs) in
one-dimension, of reaction-convection-diffusion type, whose solution
contains boundary layers (see, e.g. \cite{rst}). The approximation of SPPs has
received a lot of attention in the last few decades, mainly using finite
differences (FDs) and finite elements (FEs) on \emph{layer adapted meshes}
(see, e.g. \cite{L}). Various formulations and results are available in
the literature, both theoretical and computational \cite{L}. One method that has not,
to our knowledge, been applied to general SPPs is \emph{Isogeometric Analysis} 
(IGA). Since the introduction of IGA by T. R. Hughes et.~al. \cite{Hughes}, 
the method has been successfully applied to a large number of problem
classes. Even though much attention has been given to convection-dominated
problems \cite{HughesBook}, the method has not been applied, as far as we know, 
to a typical singularly perturbed problem, such as (\ref{de})--(\ref{bc}) ahead.

Our goal in this article is to study the appication of IGA to SPPs and in
particular the approximation of the solution to (\ref{de})--(\ref{bc}) ahead. We
use a Galerkin formulation with B-splines as basis functions and select
appropriate knot vectors, such that as the polynomial degree increases, the
error in the approximation, measured in the maximum norm, decays
exponentially, independently of any singular perturbation parameter(s).
This is the analog of performing $p$ refinement in the
FEM.

The rest of the paper is organized as follows: in Section \ref{section:iga}
we give a brief introduction to IGA, as described in \cite{HughesBook}. In
Section \ref{model} we present the model problem and its regularity. In
Section \ref{form} we give the Galerkin formulation and construct the
discrete problem. Finally, Section \ref{nr} shows the results of our
numerical computations and Section \ref{conc} gives our conclusions.

With $I\subset \mathbf{R}$ an interval with boundary $\partial I$ and
measure $\left\vert I\right\vert $, we will denote by $C^{k}(I)$ the space
of continuous functions on $I$ with continuous derivatives up to order $k$.
We will use the usual Sobolev spaces $W^{k,m}(I)$ of functions on $\Omega $
with $0,1,2,...,k$ generalized derivatives in $L^{m}\left( I\right) $,
equipped with the norm and seminorm $\left\Vert \cdot \right\Vert _{k,m,I}$
and $\left\vert \cdot \right\vert _{k,m,I}\,$, respectively. When $m=2$, we
will write $H^{k}\left( I\right) $ instead of $W^{k,2}\left( I\right) $, and
for the norm and seminorm, we will write $\left\Vert \cdot \right\Vert
_{k,I} $ and $\left\vert \cdot \right\vert _{k,I}\,$, respectively. The
usual $L^{2}(I)$ inner product will be denoted by $\left\langle \cdot ,\cdot
\right\rangle _{I}$, with the subscript omitted when there is no confusion.
We will also use the space
\[
H_{0}^{1}\left( I\right) =\left\{ u\in H^{1}\left( I\right) :\left.
u\right\vert _{\partial \Omega }=0\right\} .
\]
The norm of the space $L^{\infty }(I)$ of essentially bounded functions is
denoted by $\Vert \cdot \Vert _{\infty ,I}$. Finally, the letters $C,c$ will
be used for generic positive constants,
independent of any discretization or singular perturbation parameters.

\section{Isogeometric analysis}

\label{section:iga} 
\begin{figure}[h]
\begin{center}
\includegraphics[width=0.6\textwidth]{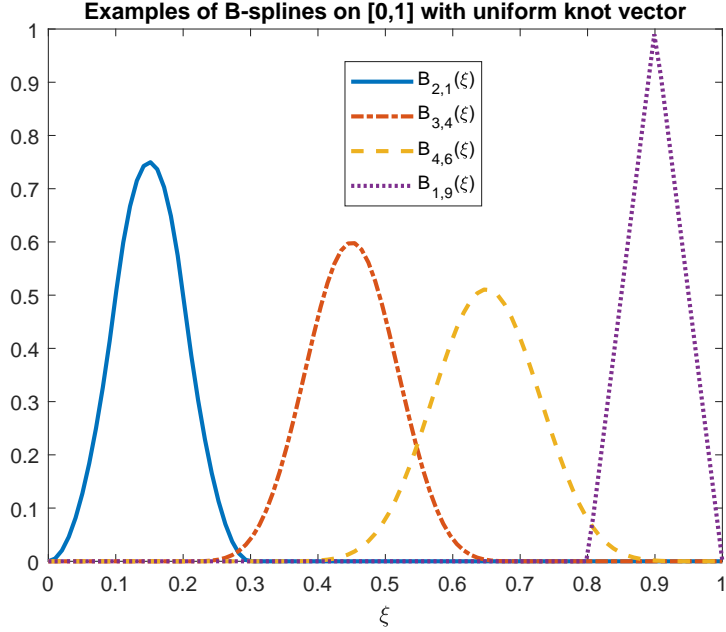}
\end{center}
\caption{Examples of B-splines using the knot vector $\Xi =[0,0.1,0.2...1]$.}
\label{F1}
\end{figure}

In this study we use B-splines as basis functions 
and follow \cite{HughesBook} closely.
To this end let $\Xi =\left\{ \xi _{1},\xi _{2},...,\xi _{n+p+1}\right\} $
be a \textbf{knot vector}, where $\xi _{i}\in \mathbf{R}$ $\ $is the $i^{th}$
knot, $i=1,2,...,N+p+1$, $p$ is the polynomial order and $N$ is the number
of basis functions used to construct the B-spline. The numbers in $\Xi $ are
non-decreasing and may be repeated, in which case we are talking about a%
\emph{\ non-uniform} knot vector. If the first and last knot values appear $%
p+1$ times, the knot vector is called \emph{open} (see \cite{HughesBook} for
more details). With a knot vector $\Xi $ in hand, the B-spline basis
functions are defined recursively, starting with piecewise constants ($p=0$): 
\[
B_{i,0}(\xi )=\left\{ 
\begin{array}{ccc}
1 & , & \xi _{i}\leq \xi <\xi _{i+1} \\ 
0 & , & otherwise%
\end{array}%
\right. .
\]
For $p=1,2,...$, they are defined by the \emph{Cox-de Boor recursion formula}
\cite{Cox,deBoor} 
\[
B_{i,p}(\xi )=\frac{\xi -\xi _{i}}{\xi _{i+p}-\xi _{i}}B_{i,p-1}(\xi )+\frac{%
\xi -\xi _{i}}{\xi _{i+p}-\xi _{i}}B_{i,p-1}(\xi ).
\]
In Figure \ref{F1} we show some of the above B-splines, obtained with the
uniform knot vector $\Xi =[0,0.1,0.2...1]$ for various polynomial degrees $p$%
.

We also mention the recursive formula for obtaining the derivative of a
B-spline \cite{HughesBook}:%
\[
\frac{d}{d\xi }B_{i,p}(\xi )=\frac{p}{\xi _{i+p}-\xi _{i}}B_{i,p-1}(\xi )-%
\frac{p}{\xi _{i+p+1}-\xi _{i+1}}B_{i+1,p-1}(\xi ).
\]

We will be considering open knot vectors, having possibly repeated entries
(other than the endpoints). If we assume we have $\xi _{1},...,\xi _{m}$
distinct knots, each having multiplicity $r_{i}$, then 
\[
\Xi =[\underset{r_{1\text{ times}}}{\underbrace{\xi _{1},...,\xi _{1}}},%
\underset{r_{2\text{ times}}}{\underbrace{\xi _{2},...,\xi _{2}}},...,%
\underset{r_{m\text{ times}}}{\underbrace{\xi _{m},...,\xi _{m}}}]
\]
and there holds $\sum_{i=1}^{m} r_{i}=N+p+1$. Since we are using open
knots, we have $r_{1}=r_{m}=p+1$. The regularity of the B-spline at each
knot $\xi _{i}$ is determined by $r_{i}$, in that the B-spline has $p-r_{i}$
continuous derivatives at $\xi _{i}$. For this reason, we define $%
k_{i}=p-r_{i}+1$ as a measure of the regularity at the knot $\xi _{i}$ and
set $\boldsymbol{k}=[k_{1},...,k_{m}].$ Note that $k_{1}=k_{m}=0$ due to the
fact we are using an open knot vector.

B-splines form a partition of unity and they span the space of piecewise
polynomials of degree $p$ on the subdivision $\left\{ \xi _{1},...,\xi
_{m}\right\} $. Each basis function is positive and has support in $[\xi
_{i},\xi _{i+p+1}]$. In the sections that follow, we will approximate the
solution to the BVP\ under consideration, using the space%
\begin{equation}
S_{\boldsymbol{k}}^{p}=span\left\{ B_{k,p}\right\} _{k=1}^{N},  \label{space}
\end{equation}
with dimension
\begin{equation}
\mathcal{N}=\dim \left( S_{\boldsymbol{k}}^{p}\right)
=m p-\sum_{i=1}^{m}k_{i}.  \label{dim}
\end{equation}%
We point out that we are using a uniform polynomial degree $p$, while we
allow for the regularity at each knot to (possibly) vary. A more general
approach would be to allow $p$ to vary as well. We will refer to $\mathcal{N}
$ as the number of degrees of freedom, DOF.


\section{The model problem \label{model}}

We will apply isogeometric analysis to the following model SPP: Find $u$
such that
\begin{eqnarray}
-\varepsilon _{1}u^{\prime \prime }(x)+\varepsilon _{2}b(x)u^{\prime
}(x)+c(x)u(x) &=&f(x)\ \text{in }I=\left( 0,1\right) ,  \label{de} \\
u(0)=u(1) &=&0\text{ },  \label{bc}
\end{eqnarray}
where $0<\varepsilon _{1},\varepsilon _{2}\leq 1$ are given parameters that
can approach zero and the functions $b,c,f$ are given and sufficiently
smooth. We assume that there exist constants $\beta ,\gamma ,\rho $,
independent of $\varepsilon _{1},\varepsilon _{2},$ such that 
$\forall \;x\in \overline{I}$ 
\begin{equation}
\label{data}
b(x)\geq \beta \geq 0\;,\;c(x)\geq \gamma >0\;,\;c(x)-\frac{\varepsilon _{2}%
}{2}b^{\prime }(x)\geq \rho >0.
\end{equation}
The structure of the solution to (\ref{de}) depends on the roots of the
characteristic equation associated with the differential operator. For this
reason, we let $\lambda _{0}(x),\lambda _{1}(x)$ be the solutions of the
characteristic equation and set
\[
\mu _{0}=-\underset{x\in \lbrack 0,1]}{\max }\lambda _{0}(x)\;,\;\mu _{1}=%
\underset{x\in \lbrack 0,1]}{\min }\lambda _{1}(x),  
\]%
or equivalently,
\begin{equation}
\label{mu}
\mu _{0,1}=\underset{x\in \lbrack 0,1]}{\min }\frac{\mp \varepsilon _{2}b(x)+
\sqrt{\varepsilon _{2}^{2}b^{2}(x)+4\varepsilon _{1}c(x)}}{2\varepsilon _{1}}.
\end{equation}
The values of $\mu _{0},\mu _{1}$ determine the strength of the boundary
layers and since $\left\vert \lambda _{0}(x)\right\vert <\left\vert \lambda
_{1}(x)\right\vert $ the layer at $x=1$ is stronger than the layer at $x=0$.
Essentially, there are three regimes \cite{L}, as shown in Table 1.

\begin{table}
\centering
\caption{Different regimes based on the relationship between $\protect%
\varepsilon_1$ and $\varepsilon_2$.}
\label{table1}       
%
%
\begin{tabular}{llll}
\hline\noalign{\smallskip}
 &  & $\mu _{0}$ & $\mu _{1}$ \\
\noalign{\smallskip}\hline\noalign{\smallskip}
convection-diffusion & $\varepsilon_{1}<<\varepsilon_{2}=1$ & $1$ & $\varepsilon_{1}^{-1}$\\
convection-reaction-diffusion $\:$ & $\varepsilon_{1}<<\varepsilon_{2}^{2}<<1 \:$ & $\varepsilon_{2}^{-1} \:$ 
& $\varepsilon_{2}/\varepsilon_{1}$\\
reaction-diffusion & $1>>\varepsilon_{1}>>\varepsilon_{2}^{2}$ & $\varepsilon_{1}^{-1/2}$ & $\varepsilon_{1}^{-1/2}$ \\
\noalign{\smallskip}\hline
\end{tabular}
\end{table}

\begin{figure}[h]
\begin{center}
\includegraphics[width=0.6\textwidth]{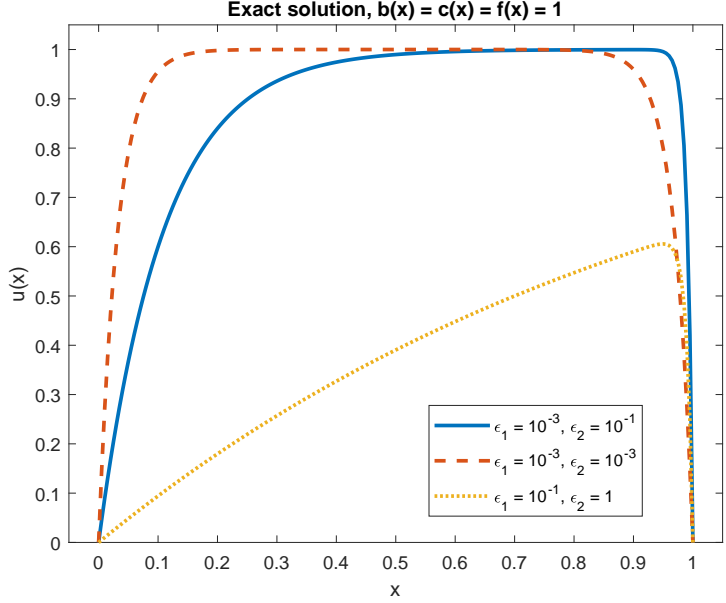}
\end{center}
\caption{Exact solution for different values of $\protect\varepsilon _{1},%
\protect\varepsilon _{2}$.}
\label{F2}
\end{figure}

We assume that $b,c,f$ are analytic functions
satisfying, for some positive constants $\gamma _{f},\gamma _{c},\gamma _{b}$ 
independent of $\varepsilon _{1},\varepsilon _{2}$, and $\forall \;n = 0, 1, 2, ... $
\[
\left\Vert f^{(n)}\right\Vert _{\infty ,I}\leq C n!\gamma
_{f}^{n}\;,\;\left\Vert c^{(n)}\right\Vert _{\infty ,I}\leq C n!\gamma
_{c}^{n}\;,\;\left\Vert b^{(n)}\right\Vert _{\infty ,I}\leq C n!\gamma
_{b}^{n}.
\]
Then, it was shown in \cite{Irene} that the solution $u$ to (\ref{de}), 
(\ref{bc}) can be decomposed into a smooth part $u_{S}$, a boundary layer part at
the left endpoint $u_{BL}^{-}$, a boundary layer part at the right
endpoint $u_{BL}^{+}$, and a (negligible) remainder, viz. 
\[
u=u_{S}+u_{BL}^{-}+u_{BL}^{+}+u_R,
\]
with
\[
\left\vert u_{S}^{(n)}(x)\right\vert \leq C K^n n!,
\]
\[
\left\vert \left(
u_{BL}^{-}\right) ^{(n)}(x)\right\vert \leq C K_1^n \mu _{0}^{n}e^{-\ell \mu
_{0}x}\;,\;\left\vert \left( u_{BL}^{+}\right) ^{(n)}(x)\right\vert \leq
C K_2^n \mu _{1}^{n}e^{-\ell \mu _{1}(1-x)},  
\]
\[
\left\Vert u_{R}\right\Vert _{\infty ,\partial I}+\left\Vert
u_{R}\right\Vert _{0,I}+\varepsilon _{1}\left\Vert u_{R}^{\prime
}\right\Vert _{0,I}\leq C \max \{e^{-\delta \varepsilon _{2}/\varepsilon
_{1}},e^{-\delta /\varepsilon _{2}}\}  ,
\]
for all $x\in \overline{I}$, where the constants $C, K, K_1, K_2, \delta >0$ depend 
\emph{only} on the data.
Figure \ref{F2} shows the behavior of the solution to (\ref{de})--(\ref{bc}%
), in all three regimes. 

If $\varepsilon_1, \varepsilon_2$ are not small (see Section \ref{nr}
for precise conditions), then no boundary layers are present and approximating
$u$ may be done using a fixed mesh (of say one element) and increasing $p$. (For 
IGA, the knot vector could simply be 
$\Xi =[0,...,0,1,...,1]$.)
If, on the other hand, $\varepsilon_1, \varepsilon_2$ are small then classical
techniques fail and the mesh must be chosen carefully.  The challenge lies 
in approximating the typical boundary layer function $\exp(-x/\varepsilon )$. 
In the context of FDs and FEs, the mesh points \textbf{must} depend on 
$\varepsilon $, as is well documented in the literature under the name 
\emph{layer-adapted} meshes \cite{L}. We expect something similar to hold
for IGA, in the sense that the knot vector \emph{must} depend on 
$\varepsilon $. We will illustrate this in Section \ref{nr}.

\section{The Galerkin formulation and the discrete problem \label{form}}

Isogeometric analysis may be combined with a number of formulations; we
choose to use Galerkin's approach, i.e. we multiply (\ref{de}) by a suitable
test function, integrate by parts and use the boundary conditions (\ref{bc}%
). \ The resulting variational formulation reads: Find $u\in H_{0}^{1}\left(
I\right) $ such that 
\begin{equation}
{\mathcal{B}}\left( u,v\right) ={\mathcal{F}}\left( v\right) \;\;\forall
\;v\in H_{0}^{1}\left( I\right) ,  \label{BuvFv}
\end{equation}%
where 
\begin{equation}
\label{Buv}
{\mathcal{B}}\left( u,v\right)  =\varepsilon _{1}\left\langle u^{\prime
},v^{\prime }\right\rangle _{I}+\varepsilon _{2}\left\langle bu^{\prime
},v\right\rangle _{I}+\left\langle cu,v\right\rangle _{I}, \:\:
{\mathcal{F}}\left( v\right)  =\left\langle f,v\right\rangle _{I}.
\end{equation}
The bilinear form ${\mathcal{B}}\left( \cdot ,\cdot \right) $ given by (\ref%
{Buv}) is \emph{coercive} (due to (\ref{data})) with respect to the \emph{%
energy norm} 
\[
\left\Vert u\right\Vert _{E,I}^{2}:=\varepsilon _{1}\left\vert u\right\vert
_{1,I}^{2}+\left\Vert u\right\Vert _{0,I}^{2},
\]%
i.e., 
\begin{equation}
{\mathcal{B}}\left( u,u\right) \geq \left\Vert u\right\Vert
_{E,I}^{2}\;\;\forall \;u\in H_{0}^{1}\left( I\right) .  \label{coercivity}
\end{equation}%
Next, we restrict our attention to a finite dimensional subspace $S\subset
H_{0}^{1}\left( I\right) $, that will be selected shortly, and obtain the
discrete version of (\ref{BuvFv}) as: find $u_{N}\in S$ such that 
\begin{equation}
{\mathcal{B}}\left( u_{N},v\right) ={\mathcal{F}}\left( v\right) \;\;\forall
\;v\in S.  \label{discrete}
\end{equation}%
The space $S$ is chosen as $S=S_{\boldsymbol{k}}^{p}$ , given by (\ref{space}%
). \ Thus, we may write the approximate solution as 
\[
u_{N}=\sum_{k=0}^{N}\alpha _{k}B_{k,p},
\]%
with $\overrightarrow{\alpha }=[\alpha _{1},...,\alpha _{N}]^{T}$ unknown
coefficients, and subsitute in (\ref{discrete}) to obtain the linear system
of equations%
\begin{equation}
\underset{M\in \mathbf{R}^{N\times N}}{\underbrace{\left( \varepsilon
_{1}A_{1}+\varepsilon _{2}A_{2}+A_{0}\right) }}\overrightarrow{\alpha }=%
\overrightarrow{F},  \label{system}
\end{equation}%
where%
\[
\lbrack A_{1}]_{i,j}=\int_{I}B_{i,p}^{\prime }(\xi )B_{j,p}^{\prime }(\xi
)d\xi \; ,\; 
\lbrack A_{2}]_{i,j}=\int_{I}B_{i,p}^{\prime }(\xi )B_{j,p}(\xi )d\xi ,
\]
\[
\lbrack A_{0}]_{i,j}=\int_{I}B_{i,p}(\xi )B_{j,p}(\xi )d\xi \; , \; 
\lbrack \overrightarrow{F}]_{i=}\int_{I}B_{i,p}(\xi )f(\xi )d\xi ,
\]
for $i,j=1,...,N.$ The linear system (\ref{system}) has a unique solution, due
to the fact that the coefficient matrix $M$ in (\ref{system}) is
non-singular. To see this, let $0\neq v\in S_{\boldsymbol{k}}^{p}$ be
arbitrary and write it as $v=\sum_{k=0}^{N}\beta _{k}B_{k,p}$, with the
coefficients $\beta _{k}$ not all zero. From (\ref{coercivity}), we have%
\begin{eqnarray*}
0 &<&\left\Vert \sum_{k=0}^{N}\beta _{k}B_{k,p}\right\Vert _{E,I}^{2}\;\leq {%
\mathcal{B}}\left( \sum_{k=0}^{N}\beta _{k}B_{k,p},\sum_{\ell =0}^{N}\beta
_{\ell }B_{\ell ,p}\right)  \\
&\leq &\sum_{k=0}^{N}\sum_{\ell =0}^{N}\beta _{k}{\mathcal{B}}\left(
B_{k,p},B_{\ell ,p}\right) \beta _{\ell }=\overrightarrow{\beta }^{T}M%
\overrightarrow{\beta },
\end{eqnarray*}%
which shows that $M$ is positive definite, hence invertible.

We close this section by mentioning that in our implementation of the
method, the entries in the matrices in (\ref{system}), i.e. integrals of
B-splines, are computed numerically to any desired accuracy (using MATLAB's 
\texttt{integrate} command).

\section{Numerical results}

\label{nr}

In this section we present the results of numerical computations for three
examples with known exact solution -- this makes our results reliable. We
will `mimic' the FEM recommendations for such problems (see, e.g., \cite%
{Irene} and the references therein), and select our open knot vector (for
the interval $I=(0,1)$) as follows:

With $\mu _{0},\mu _{1}$ given by (\ref{mu}), if $p \mu_1^{-1} \ge 1/2$ then 
\begin{equation}
\label{unif_knot}
\Xi =[\underset{p+1\text{ times}}{\underbrace{0,...,0}}, 
\underset{p+1\text{ times}}{\underbrace{1,...,1}}\}].
\end{equation}
If $p \mu_0^{-1} < 1/2$ then, for \textbf{reaction-diffusion}
\begin{equation}
\label{knot_rd}
\Xi =[\underset{p+1\text{ times}}{\underbrace{0,...,0}}, p_{\max
}\varepsilon^{1/2}_{1},1-p_{\max }\varepsilon^{1/2}_{1},\underset{p+1\text{
times}}{\underbrace{1,...,1}}\}],
\end{equation}
for \textbf{convection-diffusion}
\begin{equation}
\label{knot_cd}
\Xi =[\underset{p+1\text{ times}}{\underbrace{0,...,0}},1-p_{\max
}\varepsilon _{1},\underset{p+1\text{ times}}{\underbrace{1,...,1}}\}],
\end{equation}
and for \textbf{reaction-convection-diffusion}
\begin{equation}
\label{knot_rcd}
\Xi =[\underset{p+1\text{ times}}{\underbrace{0,...,0}},p_{\max }\mu
_{0}^{-1},1-p_{\max }\mu _{1}^{-1},\underset{p+1\text{ times}}{\underbrace{%
1,...,1}}\}],
\end{equation}
where $p=1,2,...,p_{\max }$
is the polynomial degree, which we change to improve accuracy. 
The number of degrees of freedom in each case is given by $DOF =
3p$ and we take $p_{max} = 10$ for the computations.

We will be measuring the percentage relative error in the maximum norm,
\[
Error=100\times \frac{\left\Vert u-u_{N}\right\Vert _{\infty ,I}}{\left\Vert
u\right\Vert _{\infty ,I}},
\]
which we will estimate as follows:
\[
Error\approx 100\times \underset{k=1,...,K}{\max }\left\vert
u(x_{k})-u_{N}(x_{k})\right\vert /\underset{k=1,...,K}{\max }\left\vert
u(x_{k})\right\vert ,
\]
where $\left\{ x_{k}\right\} _{k=1}^{K}\in I$ are points in $(0,1),$ chosen
uniformly in the layer region and outside -- we use $K=400$ in each region
for our experiments below. We choose to use the maximum norm as an error
measure, because the energy norm is `not balanced' for reaction-diffusion
problems (see \cite{MX} and the references therein).

The examples that follow cover all three regimes, and try to answer the
question of how the method performs as $\varepsilon _{1},\varepsilon
_{2}\rightarrow 0$. To reduce the error, we increase the dimension of the
space by increasing $p$, hence strictly speaking, we are performing $p$%
-refinement. (In the FEM literature this has been referred to as $hp$%
-refinment \cite{ss}.)

\vspace{0.5cm}

\begin{figure}[h]
\begin{center}
\includegraphics[width=0.7\textwidth]{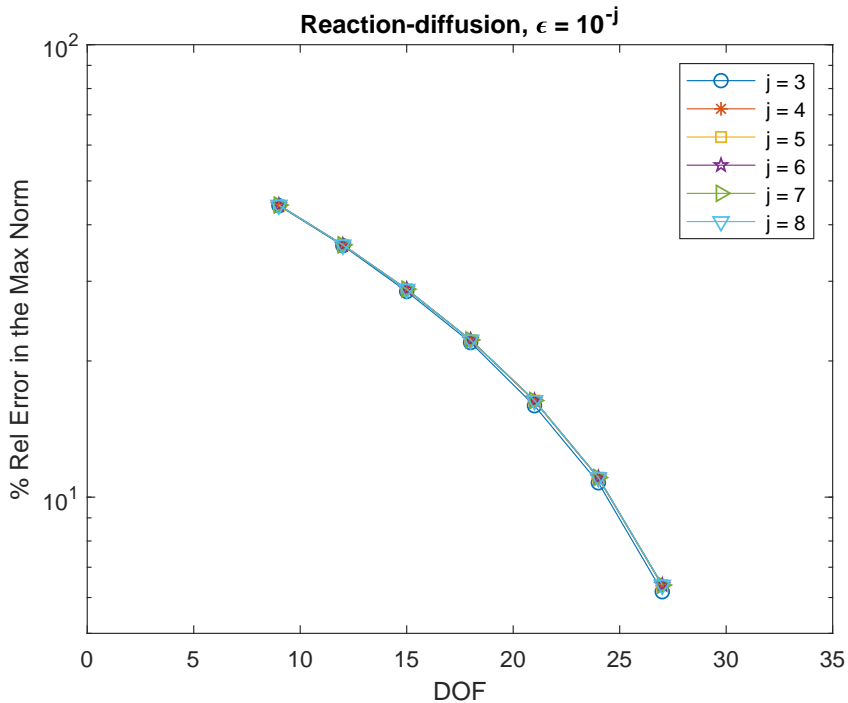}
\end{center}
\caption{Maximum norm convergence for Example 1, using the knot vector (\ref{knot_rd}).}
\label{F4}
\end{figure}

\begin{figure}[h]
\begin{center}
\includegraphics[width=0.7\textwidth]{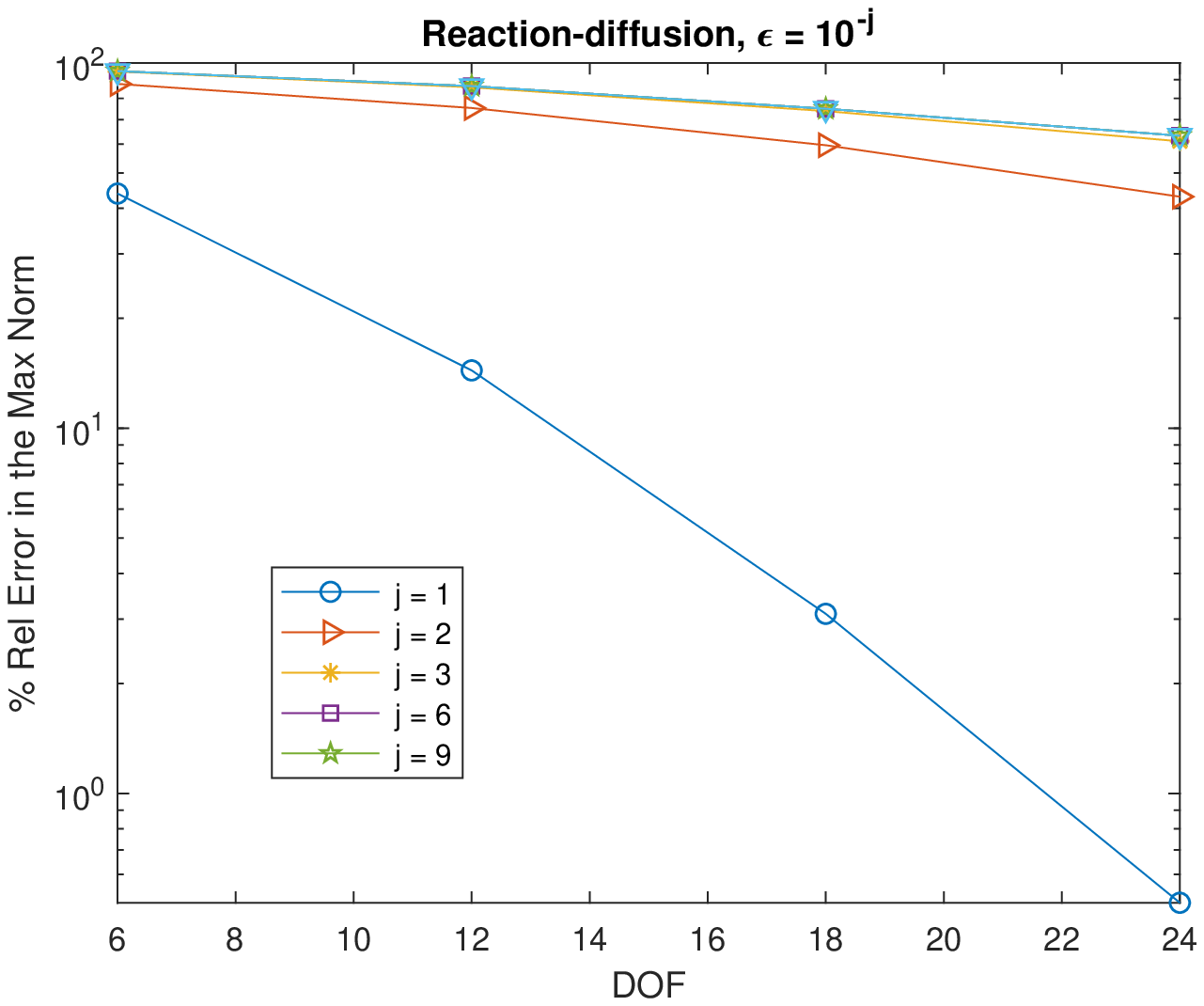}
\end{center}
\caption{Maximum norm convergence for Example 1, using the knot vector (\ref{unif_knot}).}
\label{F4b}
\end{figure}

\textbf{Example 1:} We consider (\ref{de}), (\ref{bc}) with $%
b(x)=0,c(x)=f(x)=1$, which makes the problem reaction-diffusion with $%
\mu_{0}=\mu _{1}=\varepsilon _{1}^{-1/2}$.  Figure \ref{F4} shows the
percentage relative error measured in the maximum norm, versus the number of
degrees of freedom $DOF$ (cf. (\ref{dim})) in a semi-log scale, and Table %
\ref{table2} lists the errors. The fact that we see straight lines indicates
the exponential convergence of the method, while the robustness is verified
since the straight lines coincide. We also show, in Figure \ref{F4b}, the
results of using a mesh that does \emph{not} depend on $\varepsilon_1$;
in particular we use (\ref{unif_knot}), $\varepsilon_1 = 10^{-j}, j=1,2,3,6,9,12$ and  
we increase $p$. As can be seen from the figure, for large $\varepsilon_1$ the 
method yields good results, but as $\varepsilon_1 \rightarrow 0$, the results
deteriorate and we (basically) have no convegence.

\begin{table}
\centering
\caption{Percentage relative error in the maximum norm for Example 1, using the knot vector (\ref{knot_rd}).}
\label{table2}       
%
%
\begin{tabular}{lllllll}
\hline\noalign{\smallskip}
\hspace{0.5cm}$\varepsilon_1:$ & $10^{-6}$ & $10^{-8}$ & $10^{-10}$ & $10^{-12}$ & $10^{-14}$ & $10^{-16}$ \\
$DOF$ \\
\noalign{\smallskip}\hline\noalign{\smallskip}
6 & $44.05$ & $44.19$ & $44.21$ & $44.22$ & $44.21$ & $44.21$\\
9 & $35.96$ & $36.09$ & $36.11$ & $36.12$ & $36.11$ & $36.11$\\
12 & $28.48$ & $28.78$ & $28.81$ & $28.82$ & $28.82$ & $28.82$\\
$15$ & $21.96$ & $22.23$ & $22.25$ & $22.26$ & $22.26$ & $22.26$ \\ 
$18$ & $15.92$ & $16.31$ & $16.34$ & $16.35$ & $16.35$ & $16.35$ \\ 
$21$ & $10.75$ & $11.01$ & $11.04$ & $11.04$ & $11.04$ & $11.04$ \\ 
$24$ & $6.18$ & $6.37$ & $6.39$ & $6.39$ & $6.39$ & $6.39$ \\ 
$27$ & $2.42$ & $2.41$ & $2.42$ & $2.41$ & $2.41$ & $2.41$ \\
\noalign{\smallskip}\hline
\end{tabular}
\end{table}

\begin{figure}[h]
\begin{center}
\includegraphics[width=0.7\textwidth]{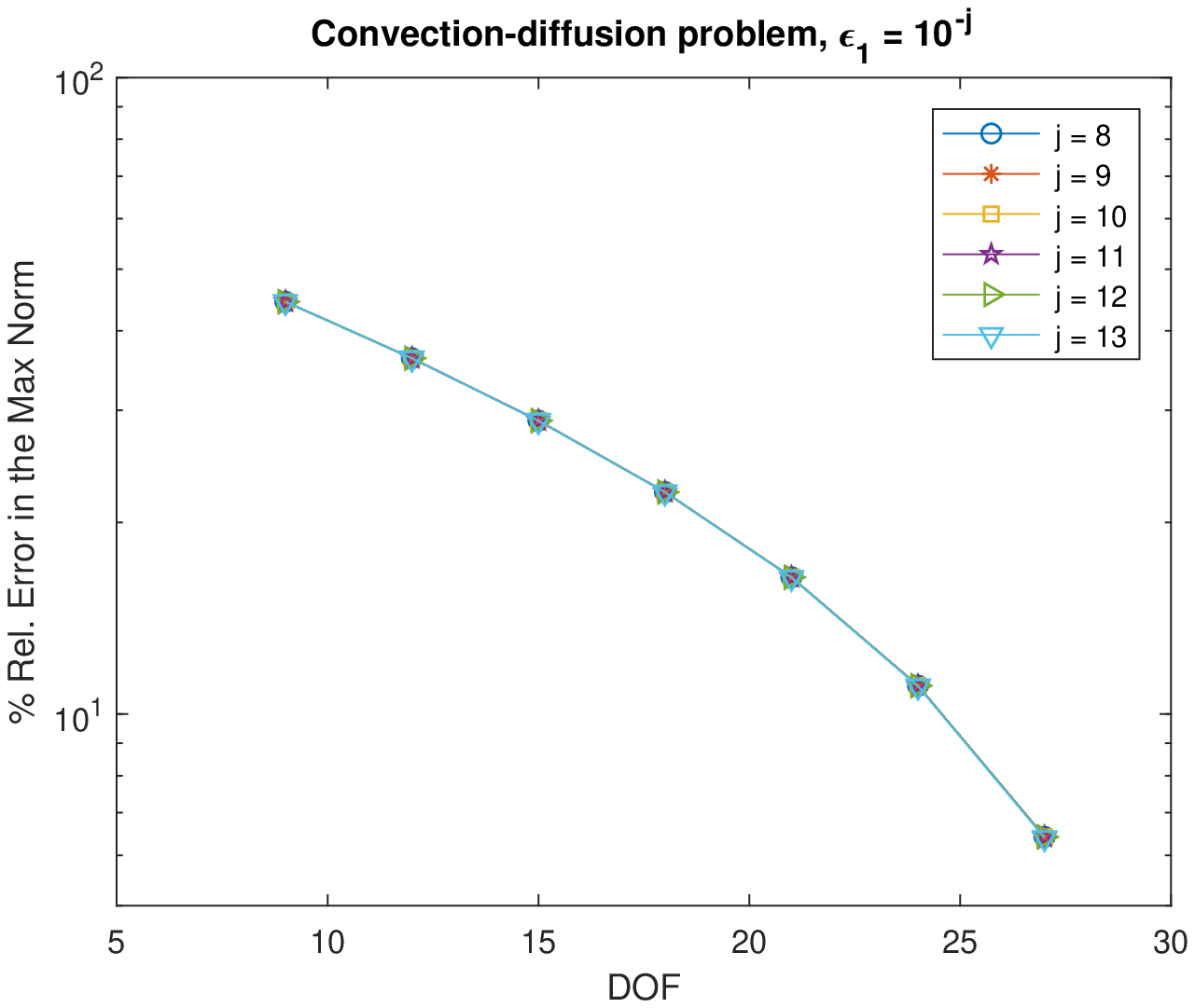}
\end{center}
\caption{Maximum norm convergence for Example 2, using the knot vector (\ref{knot_cd}).}
\label{F5}
\end{figure}

\textbf{Example 2: }We next consider (\ref{de}), (\ref{bc}) with $%
b(x)=c(x)=f(x)=\varepsilon _{2}=1$, which makes the problem
convection-diffusion with $\mu _{0}=1,\mu _{1}=\varepsilon _{1}^{-1}.$ In
Figure \ref{F5} we show the error in the maximum norm versus the number of
degrees of freedom, in a semi-log scale, for different values of $%
\varepsilon _{1}$. The errors are listed in Table \ref{table3}. Once again
we observe robust exponential convergence. The same result as Example 1
is obtained when we use the knot vector (\ref{unif_knot}), hence this is not shown
here.

\begin{table}
\centering
\caption{Percentage relative error in the maximum norm for Example 2, using the knot vector (\ref{knot_cd}).}
\label{table3}       

\begin{tabular}{lllllll}
\hline\noalign{\smallskip}
\hspace{0.5cm}$\varepsilon_1:$& $10^{-8}$ & $10^{-9}$ & $10^{-10}$ & $10^{-11}$ & $10^{-12}$ & $10^{-13}$ \\
$DOF$ \\
\noalign{\smallskip}\hline\noalign{\smallskip}
$6$ & $44.44$ & $44.45$ & $44.45$ & $44.44$ & $44.44$ & $44.44$ \\
$9$ & $36.23$ & $36.22$ & $36.22$ & $36.22$ & $36.22$ & $36.21$ \\
$12$ & $28.91$ & $28.90$ & $28.91$ & $28.91$ & $28.90$ & $28.89$ \\
$15$ & $22.32$ & $22.32$ & $22.33$ & $22.32$ & $22.32$ & $22.30$ \\
$18$ & $16.40$ & $16.41$ & $16.40$ & $16.41$ & $16.40$ & $16.38$ \\  
$21$ & $11.08$ & $11.09$ & $11.08$ & $11.08$ & $11.08$ & $11.05$ \\
$24$ & $6.41$ & $6.41$ & $6.41$ & $6.42$ & $6.41$ & $6.38$ \\  
$27$ & $2.42$ & $2.42$ & $2.42$ & $2.42$ & $2.42$ & $2.40$ \\
$30$ & $1.92$ & $1.92$ & $1.92$ & $1.93$ & $1.93$ & $1.95$  \\
\noalign{\smallskip}\hline
\end{tabular}
\end{table}

\begin{figure}[h]
\begin{center}
\includegraphics[width=0.7\textwidth]{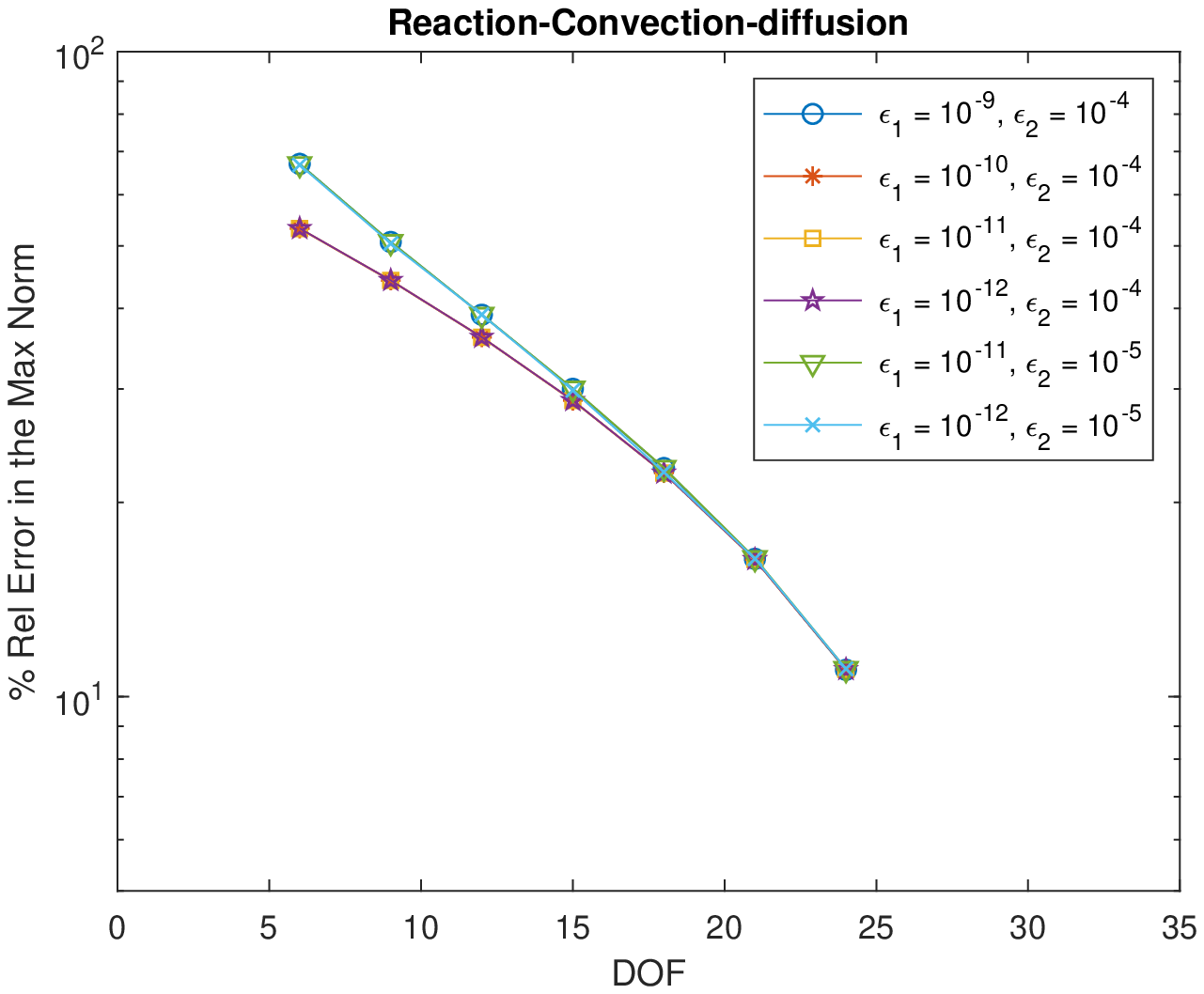}
\end{center}
\caption{Maximum norm convergence for Example 3, using the knot vector (\ref{knot_rcd}).}
\label{F6}
\end{figure}

\textbf{Example 3: }We finally consider (\ref{de}), (\ref{bc}) with $%
b(x)=c(x)=f(x)=1$, and choose the values of $\varepsilon
_{1},\varepsilon_{2} $ to satisfy $\varepsilon _{1}<<\varepsilon _{2}^{2}$,
so that the problem becomes convection-reaction-diffusion with $\mu
_{0}=\varepsilon_{2}^{-1},\mu _{1}=\varepsilon _{2}\varepsilon _{1}^{-1}$.
In Figure \ref{F6} we show the convergence of the method, measured in the
maximum norm, for different values of $\varepsilon _{1},\varepsilon _{2}$.
Table \ref{table4} shows the actual errors. We observe robust exponential
convergence in this final case as well.

\begin{table}
\centering
\caption{Percentage relative error in the maximum norm for Example 3, using the knot vector (\ref{knot_rcd}).}
\label{table4}       
\begin{tabular}{lllllll}
\hline\noalign{\smallskip}
& $\varepsilon_1=10^{-9}$ & $\varepsilon_1=10^{-10}$ & $\varepsilon_1=10^{-11}$
& $\varepsilon_1=10^{-12}$ & $\varepsilon_1=10^{-11}$ & $%
\varepsilon_1=10^{-12}$ \\ 
$DOF$ & $\varepsilon_2=10^{-4}$ & $\varepsilon_2=10^{-4}$ & $%
\varepsilon_2=10^{-4}$ & $\varepsilon_2=10^{-4}$ & $\varepsilon_2=10^{-5}$ & %
$\varepsilon_2=10^{-5}$ \\
\noalign{\smallskip}\hline\noalign{\smallskip}
$6$ & $66.98$ & $53.17$ & $53.18$ & $53.20$ & $66.97$ & $66.74$ \\ 
$9$ & $50.68$ & $44.18$ & $44.17$ & $44.19$ & $50.71$ & $50.41$ \\ 
$12$ & $39.09$ & $36.07$ & $36.08$ & $36.09$ & $39.11$ & $39.10$ \\ 
$15$ & $30.01$ & $28.77$ & $28.78$ & $28.79$ & $30.05$ & $29.84$ \\ 
$18$ & $22.61$ & $22.22$ & $22.23$ & $2.22$ & $22.65$ & $22.29$ \\ 
$21$ & $16.35$ & $16.31$ & $16.30$ & $16.32$ & $16.41$ & $16.36$ \\ 
$24$ & $11.00$ & $11.02$ & $11.03$ & $11.01$ & $11.05$ & $11.04$ \\ 
$27$ & $6.37$ & $6.38$ & $6.38$ & $6.38$ & $6.39$ & $6.39$ \\ 
$30$ & $2.43$ & $2.43$ & $2.43$ & $2.42$ & $2.41$ & $2.41$  \\
\noalign{\smallskip}\hline
\end{tabular}
\end{table}

\section{Conclusions\label{conc}}

In this article we studied the performance of IGA for one-dimensional
reaction-convection-diffusion problems with two small parameters. We
observed that if the knot vector is chosen appropriately and depending on
the singular perturbation parameter(s), then $p$-refinement yields robust,
exponential rates of convergence. The theoretical justification of what we
have observed will appear in \cite{X}.

As a next step, we intend to study the two-dimensional analogs, as well as
higher order operators. In particular, we are investigating the use of IGA
to $4^{th}$ order SPPs in two-dimensions, where their use is one of the few
available choices for obtaining an approximation, in curvilinear
two-dimensional domains.

%
%
%
%
%

%

\begin{thebibliography}{99}
%
%
%
%
%
%

\bibitem{HughesBook} J. A. Cottrell, T. R. Hughes and Y. Basilevs, \emph{%
Isogeometric Analysis: Toward integration of CAD and FEA}, Wiley and Sons 
(2009).

\bibitem{Beirao-et-al} L Beir\~{a}o da Veiga, A. Buffa, J Rivas and G.
Sangalli, \emph{Some estimates for h-p-k refinement in Isogeometric Analysis}%
, Numer. Math., \textbf{118}, 271--305 (2011).

\bibitem{Cox} M. G. Cox, \emph{The numerical evaluation of B-splines},
Technical report, National Physics Laboratory DNAC 4 (1971).

\bibitem{deBoor} C. De Boor, \emph{On calculation with B-splines}, Journal
of Approximation Theory, \textbf{6}, 50--62 (1972).

\bibitem{Hughes} T. R. Hughes, J. A. Cottrell and Y. Basilevs, \emph{%
Isogeometric analysis: CAD, finite elements, NURBS, exact geometry and mesh
refinement}, Comput. Meth. Appl. Mech. Eng., \textbf{194} (2005) 4125--4195.

\bibitem{L} T. Lin{\ss }, \emph{Layer-adapted meshes for
reaction-convection-diffusion problems}, Lecture Notes in Mathematics 1985,
Springer-Verlag, 2010.

\bibitem{MX} J. M. Melenk and C. Xenophontos, \emph{Robust exponential
convergence of hp-FEM in balanced norms for singularly perturbed
reaction-diffusion equations}, Calcolo, \textbf{53} (2016) 105--132.

\bibitem{rst} H.-G. Roos, M. Stynes, and L. Tobiska. {\emph{Robust numerical
methods for singularly perturbed differential equations}}, volume~24 of 
\emph{Springer Series in Computational Mathematics}. Springer-Verlag,
Berlin, second edition, 2008. Convection-diffusion-reaction and flow
problems.

\bibitem{ss} C. Schwab and M. Suri, \emph{The p and hp version of the finite
element method for problems with boundary layers}, Math. Comp., \textbf{65}
(1996) 1404--1429.

\bibitem{Irene} I. Sykopetritou, \emph{An hp finite element method for a
second order singularly perturbed boundary value problem with two small
parameters}, M.Sc Thesis, Department of Mathematics \& Statistics,
University of Cyprus, 2018.

\bibitem{X} C. Xenophontos and I. Sykopetritou, \emph{Isogeometric analysis for singularly
perturbed problems in 1-D: error estimates}, submitted (2019).


\end{thebibliography}
%



\end{document}